\documentclass[11pt]{amsart}
\usepackage{amssymb, latexsym}
\theoremstyle{plain}
\newtheorem{theorem}{Theorem}

\newtheorem{proposition}{Proposition}

\numberwithin{equation}{section}

\begin{document}

\title[  Non-uniform distributions on the set of parking functions ] {A family of non-uniform distributions
on the set of parking functions generated by random permutations }

\author{Ross G. Pinsky}

%\noindent  pinsky@math.technion.ac.il\ \ \ \ tel: 972-4-829-4083\ \ \  fax: 972-4-829-3388

\address{Department of Mathematics\\
Technion---Israel Institute of Technology\\
Haifa, 32000\\ Israel}
\email{ pinsky@technion.ac.il}

\urladdr{https://pinsky.net.technion.ac.il/}

\date{}

\begin{abstract}
We introduce a rather natural family of non-uniform distributions on $PF_n$, $n\in\mathbb{N}$, the set of parking functions of length $n$.
One of the motivations for this comes from a similar situation in the context of integer partitions.
For a permutation $\sigma\in S_n$ and for $j\in[n]$, let  $I_{n,<j}(\sigma)$ denote  the number of inversions in $\sigma$ that involve the number $j$ and a number less than $j$.
Let $\tilde I_{n,<j}(\sigma)=I_{n,<j}(\sigma)+1$.
The  map
$(\sigma,\tau)\to\left(\tilde I_{n,<\tau_1}(\sigma),\cdots, \tilde I_{n,< \tau_n}(\sigma)\right)$
maps $S_n\times S_n$ onto $PF_n$.
Consider the family of distributions
$P_n^{(q)}\times P_n$, $q\in(0,\infty)$, on $S_n\times S_n$, where $P_n$ is the uniform distribution on $S_n$ and
$P_n^{(q)}$ is the Mallows distribution with parameter $q$ on $S_n$. The Mallows distributions are defined by  exponential tilting via the inversion statistic.
For each $q>0$, the above map along with the distribution  $P_n^{(q)}\times P_n$ induces an exchangeable distribution $\mathcal{P}_n^{(q)}$ on $PF_n$.
We study the asymptotic behavior of two fundamental statistics of parking functions  under the family of distributions $\mathcal{P}_n^{(q)}$.
\end{abstract}

\maketitle
\section{Introduction and Statement of Results}\label{intro}

In this paper, we introduce a rather natural family of non-uniform distributions on $PF_n$, $n\in\mathbb{N}$, the set of parking functions of length $n$.
One of the motivations for this comes from a similar situation in the context of  $\Lambda_n$, $n\in\mathbb{N}$, the set of partitions of $n$, as will be described at the end of this section.
We begin by recalling the definition of a parking function.
Consider a row of  $n$ parking spaces on a one-way street. A line of $n$ cars, numbered from 1 to $n$, attempt to park, one at a time. The $i$th car's preferred space is spot number $\pi_i\in[n]$.
If this space is already taken, then car $i$ proceeds forward and parks in the first available space, if one exists. If the  car is unable to park, it exits the street.
A sequence $\pi=\{\pi_i\}_{i=1}^n$ is called a parking function on $[n]$ if all $n$ cars are able to park. It is easy to see that $\pi$ is a parking function if and only if
\begin{equation}\label{parkfuncdef}
|\{i:\pi_i\le j\}|\ge j,\ \text{for all}\ j\in[n].
\end{equation}
 It is well-known that
$|PF_n|=(n+1)^{n-1}$. There are a number of proofs of this result; a particularly  elegant one due to Pollack can be found in
\cite{FR}. Parking functions are a fundamental combinatorial object; for a survey of
 parking functions and their generalizations and some applications, see the article by Yan \cite{Y}.

One can consider random parking functions by placing the uniform probability measure on $PF_n$.
A fundamental fact about a parking function is that it remains a parking function after any permutation of its coordinates; thus, the uniform measure on $PF_n$ is exchangeable.
A study of random parking functions was initiated by Diaconis and Hicks in the fundamental paper \cite{DH}, where many interesting properties of random parking functions are proved; see also the article by Bellin \cite{Be}.
In this paper we consider a family of non-uniform, exchangeable distributions on $PF_n$ that arises naturally through the inversion statistic for permutations.

Denote by $S_n$ the set of permutations of $[n]$, and write in one-line notation $\sigma=\sigma_1\cdots\sigma_n\in S_n$. Recall that the  inversion statistic $I_n$ on $S_n$ is defined by
$$
I_n(\sigma)=\sum_{1\le i<j\le n}1_{\sigma_j<\sigma_i}=\sum_{1\le i<j\le n}1_{\sigma^{-1}_j<\sigma^{-1}_i}.
$$
(Note that $\sigma^{-1}_k$ is the position of $k$ in $\sigma$.)
 For $j\in[n]$, let  $I_{n,<j}(\sigma)$ denote  the number of inversions in $\sigma$ that involve the number $j$ and a number less than $j$; that is, $I_{n,<j}(\sigma)=\sum_{1\le  i<j}1_{\sigma^{-1}_j<\sigma^{-1}_i}$. Of course $I_{n,<1}\equiv0$, and
 $I_n=\sum_{j=1}^nI_{n,<j}$.

For convenience, define $\tilde I_{n,<j}:=I_{n,<j}+1$.
We have
\begin{equation}\label{inversionnumbers}
1\le \tilde I_{n,<j}(\sigma)\le j,\ j\in[n],\sigma\in S_n.
\end{equation}
As is well known, the map  that takes $\sigma\in S_n$ to
the vector $\left(\tilde I_{n,<1}(\sigma),\cdots, \tilde I_{n,< n}(\sigma)\right)\in [1]\times\cdots\times[n]$ is a bijection.
(The vector $\left(\tilde I_{n,<1}(\sigma),\cdots, \tilde I_{n,< n}(\sigma)\right)$ or an equivalent variant is called  the Lehmer code of the permutation.)

In light of \eqref{inversionnumbers} and \eqref{parkfuncdef}, it follows that
\begin{equation}\label{bijection}
\left(\tilde I_{n,<\tau_1}(\sigma),\cdots, \tilde I_{n,< \tau_n}(\sigma)\right)\in PF_n,
\text{for each pair}\ (\sigma,\tau)\in S_n\times S_n.
\end{equation}
  Let $P_n$ denote the uniform probability measure on $S_n$, and let $Q_n$ be an arbitrary probability measure on $S_n$. Then  \eqref{bijection} and  the probability measure
 $Q_n\times P_n$  on $S_n\times S_n$ induce an exchangeable probability measure on $PF_n$.
We will consider the probability measures $P_n^{(q)}\times P_n$, $q\in(0,\infty)$, on $S_n\times S_n$, where
$P_n^{(q)}$ is the Mallows measure with parameter $q$. The Mallows measures are defined by  exponential tilting via the inversion statistic.
More precisely,
\begin{equation*}\label{Mallows}
P_n^{(q)}(\sigma)=\frac{q^{I_n(\sigma)}}{K_n},\ \sigma\in S_n,\ q\in(0,\infty),
\end{equation*}
where the normalization constant $K_n$ is given by $K_n=\frac{(1-q)^n}{\prod_{j=1}^n(1-q^j)}$.
Of course, $P_n^{(1)}$ is the uniform distribution on $S_n$. For $q>1$, the probability $P_n^{(q)}$ favors permutations with many inversions and for $q\in(0,1)$ it favors permutations with few inversions.
A fundamental fact about Mallows distributions is that under $P_n^{(q)}$ the random variables $\{I_{n,<j}\}_{j=1}^n$ are independent with truncated geometric distributions starting at zero \cite{P22}, and thus $\{\tilde I_{n,<j}\}_{j=1}^n$ are independent with truncated geometric distributions starting at one.
More specifically, the distributions are given by
\begin{equation}\label{MallowsI<k}
P_n^{(q)}(\tilde I_{n,<j}=i)=\begin{cases}\frac{1-q}{1-q^j}q^{i-1},\ i=1,\cdots, j, \  q\in(0,\infty)-\{1\};\\ \frac1{j}, \ i=1,\cdots, j,\  q=1.\end{cases}
\end{equation}

Let $\mathcal{P}_n^{(q)}$ denote the probability measure on the set of parking functions $PF_n$ that is induced by \eqref{bijection} and the measure
$P_n^{(q)}\times P_n$ on $S_n\times S_n$.
The fact that $\{\tilde I_{n,<j}\}_{j=1}^n$ are independent with truncated geometric distributions
 makes the analysis of certain statistics of parking functions under the measures $\mathcal{P}_n^{(q)}$ quite accessible.

We note that it is easy to simulate a pair of $P_n^{(q)}\times P_n$-random permutations.
Thus, it is easy to simulate a $\mathcal{P}_n^{(q)}$-random parking function in $PF_n$.
We are unaware of a natural way to simulate  a uniformly random parking function in $PF_n$.

We will consider  the behavior of two fundamental statistics of $PF_n$ that are very amenable to analysis  under the measures
$\mathcal{P}_n^{(q)}$. Actually we will allow $q$ to depend on $n$: $q=q_n$.  The reason for this is that small changes in $q_n$ when $q_n$ is near 1 can have a dramatic influence on the asymptotic behavior of the  measures
$P_n^{(q_n)}$ on permutations and in turn on the measures
$\mathcal{P}_n^{(q_n)}$ on parking functions. To see this influence on permutations, we note the following behavior of the inversion statistic under $P_n^{(q_n)}$, which follows from  \cite{P22} along with the well-known
fact that $P_n^{(q)}(\sigma)=P_n^{(\frac1q)}(\sigma^\text{rev})$, where the reversal $\sigma^\text{rev}$ of $\sigma$ is defined by $\sigma^\text{rev}=\sigma_n\sigma_{n-1}\cdots\sigma_1$.
\begin{equation}\label{inverexp}
\begin{aligned}
&E_n^{(1)}I_n=\frac{n(n-1)}4;\\
&E_n^{(1-\frac cn)}I_n\sim k_c^-n^2,\ c>0, \text{where}\ k_c^-\in(0,\frac14)\ \text{is  an explicit constant}\\
&\text{satisfying}\ \lim_{c\to0}k_c^-=\frac14;\ \lim_{c\to\infty}k_c^-=0;\\
&E_n^{(1-\frac c{n^\alpha})}I_n\sim\frac{n^{1+\alpha}}c,\ c>0, \alpha\in(0,1);\\
&E_n^{(q)}I_n\sim\frac q{1-q}n,\ 0<q<1;\\
&E_n^{(1+\frac cn)}I_n\sim k_c^+n^2,\ c>0, \text{where}\ k_c^+\in(\frac14,\frac12)\ \text{is  an explicit constant}\\
&\text{satisfying}\ \lim_{c\to0}k_c^+=\frac14;\ \lim_{c\to\infty}k_c^+=\frac12;\\
&E_n^{(1+\frac c{n^\alpha})}I_n=\frac{n^2}2-\frac{n^{1+\alpha}}c+o(n^{1+\alpha}),\ c>0, \alpha\in(0,1);\\
&E_n^{(q)}I_n=\frac {n(n-1)}2-\frac1{q-1}n+o(n),\ q>1.
\end{aligned}
\end{equation}

We now turn to our results for  the measures $\mathcal{P}_n^{(q_n)}$.
We first consider the asymptotic  behavior of the coordinate $\pi_1$ in $\pi=(\pi_1,\cdots, \pi_n)$ under $\mathcal{P}_n^{(q_n)}$ (by exchangeability, all of the coordinates $\pi_j$ have the same distribution). In the notation, we suppress the dependence of $\pi_1$ on $n$.
 For comparison, we first  state the following theorem which encapsulates some results from \cite{DH} in the case of the uniform distribution on parking functions.
 \medskip

\bf\noindent Theorem DH\rm\ (Diaconis-Hicks).\it\ Consider $\pi_1$ under  $\mathcal{P}_n^\text{\rm unif}$, the uniform measure on the set $PF_n$ of parking functions of length $n$.

\noindent i. $\frac{\pi_1}n$ converges in the total variation norm to the uniform distribution on $[0,1]$;

\noindent ii. $\mathcal{E}_n^{\text{\rm\ unif}}\pi_1=\frac n2-\frac{\sqrt{2\pi}}4n^\frac12(1+o(1))$;

\noindent iii. Let $X$ be a random variable with the Borel distribution: $P(X=j)=e^{-j}\frac{j^{j-1}}{j!},\ j=1,2,\cdots$. Then
$$
\begin{aligned}
&\mathcal{P}_n^\text{\rm unif}(\pi_1=k)\sim\frac{1+P(X\ge k)}n, \ \text{for fixed}\ k=1,2,\cdots;\\
&\mathcal{P}_n^\text{\rm unif}(\pi_1=n-k)\sim\frac{P(X\le k+1)}n, \ \text{for fixed}\ k=0,1,2,\cdots.
\end{aligned}
$$
In particular, $\mathcal{P}_n^\text{\rm unif}(\pi_1=1)\sim\frac2n$ and $\mathcal{P}_n^\text{\rm unif}(\pi_1=n)\sim\frac1{en}$.
\rm
\medskip

\bf\noindent Remark.\rm\
%The theorem and its proof  show that in the ``bulk'', $\pi_1$ behaves uniformly, but in the ``corners'' the behavior is a little different, but on the same order.
From the definition of a parking function it is intuitively clear that $\mathcal{P}_n^\text{\rm unif}(\pi_1=j)$ is decreasing in $j\in[n]$; a formal proof can be found in
\cite{DH}. The above theorem and its proof show that the probabilistic advantage of smaller numbers
disappears as $n\to\infty$, except for a remnant that can be seen in the ``corners''.

\medskip

Turning to the measure $\mathcal{P}_n^{(q)}$, we
 begin with a basic and useful fact.
\begin{proposition}\label{stochdom}
Let $n\ge2$. The random variable $\pi_1$ under $\mathcal{P}_n^{(q)}$  is strictly stochastically increasing in $q$; specifically, if $q'<q''$, then
$\mathcal{P}_n^{(q')}(\pi_1\ge k)<\mathcal{P}_n^{(q'')}(\pi_1\ge k),\ \text{for}\ k=2,\cdots, n$.
\end{proposition}

It is easy to see from the definition of $\mathcal{P}_n^{(1)}$ that
$\mathcal{P}_n^{(1)}(\pi_1=j)$ is decreasing in $j\in[n]$. And it is not hard to see intuitively from the construction that this advantage to smaller numbers is stronger than it is in
the case of the uniform measure $\mathcal{P}_n^\text{\rm unif}$ on $PF_n$.
In light of Proposition \ref{stochdom}, this advantage becomes even more pronounced  when $q_n<1$. For $q_n>1$, this advantage begins to break down.
All of this will be seen quantitatively in the theorems below.
We begin with the case  $q_n\equiv1$.
\begin{theorem}\label{p1q1}
 Consider $\pi_1$ under the measure $\mathcal{P}_n^{(1)}$.

\noindent i. $\frac{\pi_1}n$ converges in distribution to the distribution  on $[0,1]$ with  density\newline $-\log x,\ x\in(0,1]$;

\noindent ii. $\mathcal{E}_n^{(1)}\pi_1=\frac{n+3}4$;

\noindent iii.  $\mathcal{P}_n^{(1)}(\pi_1=k_n)\sim\frac{\log n}n,\ \text{for}\ 1\le k_n=o(n)$;

\noindent iv. $\mathcal{P}_n^{(1)}(\pi_1=k_n)\sim\frac{-\log d}n,\ \text{for}\ k_n\sim dn,\  d\in(0,1)$;

\noindent v. $\mathcal{P}_n^{(1)}(\pi_1=n-k_n)\sim\frac{1+k_n}{n^2},\ \text{for}\ 0\le k_n=o(n)$.
\end{theorem}
\bf\noindent Remark 1.\rm\ Unlike in the case of the uniform measure on parking functions, in the case of the measure
$\mathcal{P}_n^{(1)}$,  the advantage to smaller numbers remains in the limit since by part (i) the limiting density is decreasing rather than constant.

\noindent \bf Remark 2.\rm\ Compare parts (iii) and (v) in Theorem \ref{p1q1} to part (iii) of Theorem DH above.
\
\medskip

For most of the  rest of the regimes of $q_n$, we suffice with the weak convergence and don't analyze the exceptional behavior in the corners.
We next consider $q_n=1+\frac cn$, with $c\in\mathbb{R}-\{0\}$, the range of $q_n$ for which the expected number of inversions in a $P_n^{(q_n)}$-distributed permutation  is on the order $n^2$, as it is in the  case of a uniformly random permutation, as seen in \eqref{inverexp}.
\begin{theorem}\label{p1q11n}
Consider $\pi_1$ under the measure $\mathcal{P}_n^{(1+\frac cn)}$, with $c\in \mathbb{R}-\{0\}$.
Then $\frac{\pi_1}n$ converges in distribution to the distribution on $[0,1]$ with
distribution function
\begin{equation}\label{distfunc1n}
F_c(x)=\begin{cases}x+\frac1c(e^{cx}-1)\log\frac{1-e^{-c}}{1-e^{-cx}},\ x\in[0,1],\ \text{if}\ c>0;\\
x+\frac1{|c|}(1-e^{-|c|x})\log\frac{e^{|c|}-1}{e^{|c|x}-1},\ x\in [0,1],\ \text{if}\ c<0,\end{cases}
\end{equation}
and
monotone decreasing density function
\begin{equation}\label{density1n}
f_c(x)=\begin{cases} e^{cx}\log\frac{1-e^{-c}}{1-e^{-cx}},\ x\in(0,1],\ \text{if}\ c>0;\\
e^{-|c|x}\log\frac{e^{|c|}-1}{e^{|c|x}-1},\ x\in(0,1],\ \text{if}\ c<0.\end{cases}
\end{equation}

\end{theorem}
\bf\noindent Remark.\rm\ Note that the  advantage to smaller numbers still remains in the limit even when $q=1+\frac cn>1$, since the limiting density function is monotone decreasing.
\medskip

Now we consider  $q_n=1-\frac c{n^\alpha},\ c>0$, the case that the expected number of inversions in  a $P_n^{(q_n)}$-distributed permutation  is on  the order $n^{1+\alpha}$, as seen in \eqref{inverexp}.
The appropriate scaling order  to obtain a limit in distribution is now  $n^\alpha$.
\begin{theorem}\label{p1q1nalpha}
Consider $\pi_1$ under the measure $\mathcal{P}_n^{(1-\frac c{n^\alpha})}$, with $c>0$ and $\alpha\in(0,1)$.
Then $\frac{\pi_1}{n^\alpha}$ converges in distribution to the distribution on  $[0,\infty)$ with density $ce^{-cx}$, that is, to the exponential distribution with parameter $c$.
\end{theorem}

The above theorem shows that in the regime $q_n=1-\frac c{n^\alpha}$, $\pi_1$ essentially does not take on values on order above $n^\alpha$.
The reduction in possible values for $\pi_1$ becomes even more pronounced  when we move to the case of constant $q\in(0,1)$. Now a limiting distribution is obtained without any scaling.
\begin{theorem}\label{p1q}
Consider $\pi_1$ under the measure $\mathcal{P}_n^{(q)}$, with $q\in(0,1)$.
Then
\begin{equation}\label{p1qdist}
\lim_{n\to\infty}\mathcal{P}_n^{(q)}(\pi_1=k)=(1-q)q^{k-1},\ k\in\mathbb{N}.
\end{equation}
That is, $\pi_1$ converges in distribution to the geometric distribution with parameter $1-q$.
\end{theorem}

Now we consider $q_n>1$ bounded with $q_n-1$  on  order larger than $O(\frac1n)$. In particular, this includes the case $q_n\equiv q>1$.
By \eqref{inverexp},
for such regimes one has  $E_n^{(q_n)}I_n\sim\frac{n^2}2$.
  It turns out that for any such regime, there is weak convergence
to the uniform measure on $[0,1]$,  as is  the case under the uniform distribution $\mathcal{P}_n^{\text{unif}}$ on parking functions, as seen in Theorem DH.
\begin{theorem}\label{p1qlarge}
Consider $\pi_1$ under the measure $\mathcal{P}_n^{(q_n)}$ with $q_n$ bounded and
$0<q_n-1=\omega(\frac1n)$ (equivalently, $\frac1n=o(q_n-1)$).
Then $\frac{\pi_1}n$ converges in distribution to the uniform  distribution on $[0,1]$.

Also, for fixed $q>1$,
\begin{equation}\label{cornersq}
\begin{aligned}
&\mathcal{P}_n^{(q)}(\pi_1=k)=\frac1nq^{k-1}(q-1)\sum_{j=k}^n\frac1{q^j-1},\ k=1,2,\cdots, n;\\
&\mathcal{P}_n^{(q)}(\pi_1=n-k)\sim\frac1n\left(1-\frac1{q^{k+1}}\right),\ k=0,1,\cdots.
\end{aligned}
\end{equation}

\end{theorem}
\noindent \bf Remark.\rm\ Compare \eqref{cornersq} to part (iii) of Theorem DH above.
\medskip

We now turn to another statistic of  parking functions.
Define
\begin{equation}\label{Nk}
N^{(n)}_k(\pi)=\sum_{i=1}^n1_k(\pi_i),\ \pi\in PF_n, \ k\in[n],
\end{equation}
 the number of times that $k$ appears as an entry in $\pi$.
This statistic was studied briefly in \cite{DH}. Note that $N^{(n)}_1\ge1$.
It was shown there that
under the uniform measure  $\mathcal{P}_n^\text{\rm unif}$ on $PF_n$, the statistic
$N^{(n)}_1$ converges in distribution to $X+1$, where $X$ has the Poisson distribution with parameter one.
Note that $N^{(n)}_n$ can only take on the values 0 or 1. It was noted in \cite{DH} that it is easy to show that
$N^{(n)}_n$ converges in distribution to the Bernoulli distribution with parameter $\frac1e$.
There don't seem to be any results in the literature for the limiting distribution of $N^{(n)}_{k_n}$ under the uniform measure for any cases of $k_n$ besides  $k_n\equiv1$ and $k_n=n$.
We study here the limiting distribution of $N^{(n)}_{k_n}$ under  $\mathcal{P}_n^{(q_n)}$ for various regimes of $\{k_n\}$ and various regimes of $\{q_n\}$.

\begin{theorem}\label{Nkqnear1}
Consider $N^{(n)}_k$ as in \eqref{Nk}  under the measure $\mathcal{P}_n^{(q_n)}$ with\newline $q_n-1=O(\frac1n)$.

\noindent i. For fixed $k$, $\lim_{n\to\infty}\frac{N^{(n)}_k}{\log n}\stackrel{\text{\rm dist}}{=}1$;

\noindent ii.  For  $k_n=\theta(n^\alpha)$ with  $\alpha\in(0,1)$, $\lim_{n\to\infty}\frac{N^{(n)}_{k_n}}{\log n}\stackrel{\text{\rm dist}}{=}1-\alpha$;

\noindent iii. For $k_n\sim dn,\ d\in(0,1)$ and $q_n=1+\frac cn,\ c\in\mathbb{R}$, $N^{(n)}_{k_n}$ converges in distribution to   the Poisson distribution with parameter
$\lambda_c$ given by
$$
\lambda_c=\begin{cases} e^{cd}\log\frac{1-e^{-c}}{1-e^{-cd}},\ \text{if}\ c>0;\\ e^{cd}\log\frac{e^c-1}{e^{cd}-1},\ \text{if}\ c<0;\\ -\log d,\ \text{if}\ c=0.\end{cases}
$$

\noindent iv. For $k_n$ satisfying $n-k_n=o(n)$, $N^{(n)}_{k_n}$ converges in distribution to 0;

\end{theorem}
\bf\noindent Remark.\rm\ Note that the density in \eqref{density1n} from Theorem \ref{p1q11n} is the same as the parameter of the Poisson distribution in part (iii) of the Theorem \ref{Nkqnear1}.

\begin{theorem}\label{Nkq<1}
Consider $N^{(n)}_k$ as in \eqref{Nk}  under the measure $\mathcal{P}_n^{(q)}$ with $q<1$.

\noindent i. For fixed $k$, $\lim_{n\to\infty}\frac{N^{(n)}_k}n\stackrel{\text{\rm dist}}{=}(1-q)q^{k-1}$;

\noindent ii. If $\lim_{n\to\infty} nq^{k_n}=\infty$, then $\lim_{n\to\infty}\frac{N^{(n)}_{k_n}}{nq^{k_n}}\stackrel{\text{\rm dist}}{=}\frac{1-q}q$.

\noindent iii. If $L:=\lim_{n\to\infty} nq^{k_n}\in(0,\infty)$, then $N^{(n)}_{k_n}$ converges in distribution to the Poisson distribution with parameter $\frac{1-q}qL$;

\noindent iv. If $\lim_{n\to\infty} nq^{k_n}=0$, then $\lim_{n\to\infty}N^{(n)}_{k_n}\stackrel{\text{\rm dist}}{=}0$;

\end{theorem}

\bf\noindent Remark 1.\rm\ In light of  \eqref{p1qdist} in Theorem \ref{p1q} and the exchangeability of the measure $\mathcal{P}_n^{(q)}$, part (i) of Theorem \ref{Nkq<1} is a kind of  weak law of large numbers for the dependent random variables $\{1_{\{\pi_j=k\}}\}_{j=1}^\infty$.

\noindent \bf Remark 2.\rm\ If $k_n\le a\log n$, with $a<\frac1{|\log q|}$, then (ii) holds, while if
 $k_n\ge a\log n$, with $a>\frac1{|\log q|}$, then (iv) holds. In particular, a $\mathcal{P}_n^{(q)}$-random parking function with $q<1$ will most likely include no numbers of size larger than $a\log n$, if $a>\frac1{|\log q|}$, but it will include numbers of size at least $a\log n$ if $a<\frac1{|\log q|}$.

\newpage
\begin{theorem}\label{Nkq>1}
Consider $N^{(n)}_k$ as in \eqref{Nk}  under the measure $\mathcal{P}_n^{(q)}$ with $q>1$.

\noindent i. For $k$ fixed, $N^{(n)}_k$ converges in distribution  to the distribution of
$\sum_{j=k}^\infty Z_{j;k}$, where $\{Z_{j;k}\}_{j=k}^\infty$ are independent random variables and $Z_{j;k}$ has the Bernoulli distribution with parameter $\frac{(q-1)q^{k-1}}{q^j-1}$.
The distribution is supported on $\mathbb{Z}^+$;
\

\noindent ii. For $\{k_n\}_{n=1}^\infty$ satisfying $\lim_{n\to\infty} k_n=\infty$ and $\lim_{n\to\infty}(n-k_n)=\infty$, $N^{(n)}_{k_n}$ converges in distribution to the distribution of
$\sum_{i=0}^\infty Y_i$, where $\{Y_i\}_{i=0}^\infty$ are independent and $Y_i$ has the
Bernoulli distribution with parameter $\frac{q-1}{q^{i+1}}$. The distribution is supported on $\mathbb{Z}^+$;
\

\noindent iii. For $k$ fixed, $N^{(n)}_{n-k}$ converges in distribution to the distribution of $\sum_{i=0}^k Y_i$. The distribution is supported on $\{0,1,\cdots, k+1\}$.
\

\noindent iv. $\sum_{i=0}^\infty Y_i\stackrel{\text{dist}}{=}\lim_{k\to\infty}\sum_{i=0}^kY_i\stackrel{\text{dist}}{=}\lim_{k\to\infty}\sum_{j=k}^\infty Z_{j;k}$.
Thus,
$$
\lim_{k\to\infty}\lim_{n\to\infty}N^{(n)}_{n-k}\stackrel{\text{dist}}{=}\lim_{k\to\infty}\lim_{n\to\infty}N^{(n)}_k=\sum_{i=0}^\infty Y_i.
$$
\end{theorem}

The above results suggest the following question.

\noindent \bf Open Question.\rm\ Asymptotically as $n\to\infty$, are there  regimes of $\{q_n\}$ for which  in  total variation norm one has
$\lim_{n\to\infty}||\mathcal{P}_n^{(q_n)}-\mathcal{P}_n^{\text{unif}}||_{\text{TV}}<1$, and if so, for what regimes is this limit  minimized, and what is the value of the minimum?

\medskip

In light of Theorem DH and Theorems \ref{p1q1}-\ref{p1qlarge}, it would seem that a minimum less than one, if it exists, might occur at some fixed $q>1$. (We note that it is easy to show that
 for fixed $n$ and $q\to\infty$, $\mathcal{P}_n^{(q)}$  converges in distribution to the uniform measure on the set of permutations $S_n\subset PF_n$. Thus, a minimum less than one will certainly not occur
 for unbounded sequences of $\{q_n\}$.)

 As noted in the paragraph before Theorem \ref{Nkqnear1}, under $\mathcal{P}_n^{\text{unif}}$, $N^{(n)}_1$ converges in distribution to the distribution of
  $X+1$, where $X\sim \text{Pois}(1)$, and $N^{(n)}_n$ converges in distribution to $\text{Ber}(\frac1e)$.
 Using this with Theorem \ref{Nkq>1}, we can get lower bounds on  $\lim_{n\to\infty}||\mathcal{P}_n^{(q)}-\mathcal{P}_n^{\text{unif}}||_{\text{TV}}$, for $q>1$,
 and thus on $\inf_{q>1}\lim_{n\to\infty}||\mathcal{P}_n^{(q)}-\mathcal{P}_n^{\text{unif}}||_{\text{TV}}$.
 We have $\lim_{n\to\infty}\mathcal{P}_n^{\text{unif}}(N^{(n)}_1=1)=P(X=0)=e^{-1}$, and by part (i) of Theorem \ref{Nkq>1} with $k=1$,
 $$
 \lim_{n\to\infty}\mathcal{P}_n^{(q)}(N^{(n)}_1=1)=\prod_{j=2}^\infty\left(1-\frac{q-1}{q^j-1}\right)=\frac{q-1}q.
 $$
 And we have $\lim_{n\to\infty}\mathcal{P}_n^{\text{unif}}(N^{(n)}_n=1)=\frac1e$ and by part (iii) of   Theorem \ref{Nkq>1} with $k=0$,
 $\lim_{n\to\infty}\mathcal{P}_n^{(q)}(N^{(n)}_n=n)=\frac {q-1}q$.
 Thus,
when $q=\frac e{e-1}$, the limiting probabilities of the events $\{N^{(n)}_1=1\}$ and $\{N^{(n)}_n=1\}$
under  $\mathcal{P}_n^{(q)}$ and under $\mathcal{P}_n^{\text{unif}}$ are all the same. Therefore, this information does  not yield a nontrivial lower bound
on $\inf_{q>1}\lim_{n\to\infty}||\mathcal{P}_n^{(q)}-\mathcal{P}_n^{\text{unif}}||_{\text{TV}}$.
Now consider the event  $\{N^{(n)}_1=2\}$.
We have $\lim_{n\to\infty}\mathcal{P}_n^{\text{unif}}(N^{(n)}_1=2)=P(X=1)=e^{-1}$,
and by part (i) of Theorem \ref{Nkq>1}
with $k=2$,
$$
\begin{aligned}
&\lim_{n\to\infty}\mathcal{P}_n^{(q)}(N^{(n)}_1=2)=
\sum_{l=2}^\infty P\left( Z_{1;l}=1,Z_{1;k}=0,k\in\mathbb{N}-\{1,l\}\right)=\\
&\sum_{l=2}^\infty \frac{q-1}{q^l-q}\prod_{j=2}^\infty\left(1-\frac{q-1}{q^j-1}\right)=\frac{(q-1)^2}{q^2}\sum_{l=1}^\infty\frac1{q^l-1}.
\end{aligned}
$$
So
$$
\begin{aligned}
&\liminf_{n\to\infty}||\mathcal{P}_n^{(q)}-\mathcal{P}_n^{\text{unif}}||_{\text{TV}}\ge\\
&\lim_{n\to\infty}\max\left(|\mathcal{P}_n^{(q)}(N^{(n)}_n=1)-\mathcal{P}_n^{\text{unif}}(N^{(n)}_n=1)|,|\mathcal{P}_n^{(q)}(N^{(n)}_1=2)-\mathcal{P}_n^{\text{unif}}(N^{(n)}_1=2)|\right)=\\
&\max\left(\left|\frac {q-1}q-\frac1e\right|,\left|\frac{(q-1)^2}{q^2}\sum_{l=1}^\infty\frac1{q^l-1}
-\frac1e\right|\right).
\end{aligned}
$$
Using a graphing calculator, we found that the right hand side above attains its minimum value of approximately
$0.058$ at $q\approx 1.74$.
Thus,
$$
\begin{aligned}
&\inf_{q>1}\liminf_{n\to\infty}||\mathcal{P}_n^{(q)}-\mathcal{P}_n^{\text{unif}}||_{\text{TV}}\ge0.058.
\end{aligned}
$$
We note that, in fact, we  calculated explicitly
$$
\lim_{n\to\infty}\max\left(|\mathcal{P}_n^{(q)}(N^{(n)}_n=1)-\mathcal{P}_n^{\text{unif}}(N^{(n)}_n=1)|,|\mathcal{P}_n^{(q)}(N^{(n)}_1\in A)-\mathcal{P}_n^{\text{unif}}(N^{(n)}_1\in A)|\right)
$$
for all choices of $A\subset\{1,2,3,4\}$. Denoting   the limit by $L(A,q)$,  we found that $\inf_{q>1}L(A,q)$ was the smallest in the case of $A=\{2\}$.

It is intuitively clear, and not difficult to prove, that $N^{(n)}_k$ is stochastically decreasing in $k$ under $\mathcal{P}_n^{\text{unif}}$.
  Thus, for any $\{k_n\}_{n=1}^\infty$, the limiting distribution of $N^{(n)}_{k_n}$ under $\mathcal{P}_n^{\text{unif}}$, if it exists,  is stochastically dominated
  by the distribution of $X+1$ and stochastically dominates the distribution $\text{Ber}(\frac1e)$.
  It would be interesting to know  the limiting distribution
of $N^{(n)}_{k_n}$ under $\mathcal{P}_n^{\text{unif}}$ for various choices of $\{k_n\}$.
  Such information in conjunction with Theorems \ref{Nkqnear1} and \ref{Nkq>1} would
 allow for   lower bounds   on  $\lim_{n\to\infty}||\mathcal{P}_n^{(q_n)}-\mathcal{P}_n^{\text{unif}}||_{\text{TV}}$, when $q_n>1$, in the
spirit of the calculations made above.

\medskip

We end this section with a description of one of the
 motivations for the  construction  of the  family of distributions on $PF_n$ that  appears in this  paper.
This motivation comes from a rather parallel situation for partitions of an integer.
A partition $\lambda=(\lambda_1,\cdots, \lambda_k)$ of $n\in\mathbb{N}$ is a non-increasing  sequence of positive integers, $\{\lambda_i\}_{i=1}^k$,  for some $k\in[n]$, satisfying
$\sum_{i=1}^k\lambda_i=n$. Denote by $\Lambda_n$ the set of all partitions of $n$.
The highly nontrivial asymptotic growth rate of $|\Lambda_n|$ was first proved by Hardy and Ramanujan in 1918; one has $|\Lambda_n|\sim\frac1{4\sqrt3 n}e^{\pi\sqrt{\frac{2n}3}}$.
We are unaware of   a natural way to simulate a uniformly random partition from $\Lambda_n$.
In order to study the probabilistic behavior of uniformly random partitions, Fristedt in the beautiful paper \cite{F} imbedded the problem
in a one-parameter family of probability measures on partitions of random length, based on the  generating function for partition functions.

However, there is a family of non-uniform probability measures on $\Lambda_n$ that arises very naturally via permutations, and thus is easy to simulate.
Any permutation $\sigma\in S_n$ can be decomposed uniquely into the product of disjoint cycles (in any order since the disjoint cycles commute).
The  lengths of these disjoint cycles, listed in non-increasing order, constitute a partition in $\Lambda_n$.
(For example, the permutation $\sigma=62534187\in S_8$ has cycle decomposition $\sigma=(435)(16)(78)(2)$, yielding the partition $(3,2,2,1)\in \Lambda_8$.)
Consequently, any probability measure on $S_n$ induces a probability measure on $\Lambda_n$.
Consider the family of Ewens sampling distributions, which are defined by exponential tilting via the total-number-of-cycles statistic. More precisely,
let $D_n(\sigma)$ denote the total number of disjoint cycles in $\sigma\in S_n$. For each $\theta\in(0,\infty)$, define the probability measure
$Q_n^{(\theta)}(\sigma)=\frac{\theta^{D_n(\sigma)}}{\theta^{(n)}}$, for $\sigma\in S_n$, where the normalization constant
$\theta^{(n)}$ is the rising factorial, $\theta^{(n)}=\theta(\theta+1)\cdots(\theta+n-1)$.
Of course, $Q_n^{(1)}$ is the uniform distributions on $S_n$. For $\theta>1$, the probability $Q_n^{(\theta)}$ favors permutations with many cycles and for $\theta\in(0,1)$ it favors permutations with few cycles.
Let $\mathcal{Q}_n^{(\theta)}$ denote the measure induced on $\Lambda_n$ corresponding to the Ewens sampling measure $Q_n^{(\theta)}$ on $S_n$.
Natural permutation statistics, such as for example, the total number of cycles, the number of cycles of length $j$ and the length of the longest cycle, translate into corresponding statistics for partitions---the total number of components, the number of components of length $j$ and the size of the largest component.
The asymptotic probabilistic behavior of these permutation statistics and others under $Q_n^{(\theta)}$, and thus equivalently these and other statistics for partitions under $\mathcal{Q}_n^{(\theta)}$, can be found in  the fine book \cite{ABT}.
The measure itself,
$\mathcal{Q}_n^{(\theta)}$ on $\Lambda_n$, properly scaled, converges to a Poisson-Dirichlet distribution on the set of non-negative, non-decreasing sequences  that sum to one; see \cite{ABT} and references therein.
In contrast to the situation in this paper with the parameter $q$, there is no reason to consider the Ewens sampling distributions with parameter $\theta$ depending on $n$ because
the expected value of the total number of cycles $D_n$ is asymptotic to $\theta\log n$ under $Q_n^{(\theta)}$; thus, the order remains the same for all values of $\theta$.

We note that with respect to the open question we posed above for parking functions, in the case of partitions the total variation norm between the two measures above converges to one, for all values of $\theta$.
This can be seen, for example, by comparing Theorem 10.1 to equation (36) or by comparing Theorem 10.2 to equation (37) in \cite{P21}.

In section \ref{pi1proofs}, we prove Proposition \ref{stochdom} and Theorems \ref{p1q1}-\ref{p1qlarge}, and in section \ref{Nkproofs} we prove
Theorems \ref{Nkqnear1}-\ref{Nkq>1}.

\section{Proofs of Proposition \ref{stochdom} and Theorems \ref{p1q1}-\ref{p1qlarge}}\label{pi1proofs}
In this section we prove the results concerning $\pi_1$.
From \eqref{bijection}, \eqref{MallowsI<k} and  the definition of
$\mathcal{P}_n^{(q)}$, it follows that
\begin{equation}\label{p1dist=}
\mathcal{P}_n^{(q)}(\pi_1=k)=\begin{cases}\frac1n\sum_{j=k}^n\frac{1-q}{1-q^j}q^{k-1},\ q\in (0,\infty)-\{1\};\\
\frac1n\sum_{j=k}^n\frac1j,\ q=1,
\end{cases}
\end{equation}
and consequently,
\begin{equation}\label{p1distle}
\mathcal{P}_n^{(q)}(\pi_1\le k)=\begin{cases}\frac1n\sum_{l=1}^k\sum_{j=l}^n\frac{1-q}{1-q^j}q^{l-1},\ q\in (0,\infty)-\{1\};\\
\frac1n\sum_{l=1}^k\sum_{j=l}^n\frac1j,\ q=1.
\end{cases}
\end{equation}
For $q\neq1$, we have
$$
\begin{aligned}
&\sum_{l=1}^k\sum_{j=l}^n\frac{1-q}{1-q^j}q^{l-1}=\sum_{j=1}^k\sum_{l=1}^j\frac{1-q}{1-q^j}q^{l-1}+\sum_{j=k+1}^n\sum_{l=1}^k\frac{1-q}{1-q^j}q^{l-1}=\\
&k+(1-q^k)\sum_{j=k+1}^n\frac1{1-q^j}.
\end{aligned}
$$
Similarly, for the case $q=1$ we have
$$
\sum_{l=1}^k\sum_{j=l}^n\frac1j=k+k\sum_{j=k+1}^n\frac1j.
$$
Using this, we can rewrite \eqref{p1distle} as
\begin{equation}\label{p1distfunc}
\mathcal{P}_n^{(q)}(\pi_1\le k)=\begin{cases}\frac kn+\frac1n(1-q^k)\sum_{j=k+1}^n\frac1{1-q^j},\ q\in(0,\infty)-\{1\};\\
\frac kn+\frac kn\sum_{j=k+1}^n\frac1j,\ q=1.
\end{cases}
\end{equation}
Using \eqref{p1distfunc}, we now prove  Proposition \ref{stochdom} and Theorems \ref{p1q1}-\ref{p1qlarge}.
\medskip

\it\noindent
Proof of Proposition \ref{stochdom}.\rm\
Fix $1\le k<j\le n$, and  let
$$
G(q)=\begin{cases}\frac{1-q^k}{1-q^j},\ q\in(0,\infty)-\{1\};\\  \frac kj,\ q=1.\end{cases}
$$
Note that $G$ is continuous at $q=1$.
From \eqref{p1distfunc}, it suffices to show that
$G$ is strictly decreasing.
For $q\neq1$, we have
$G'(q)=\frac{H(q)}{(1-q^j)^2}$, where
$H(q)=-(1-q^j)kq^{k-1}+(1-q^k)jq^{j-1}$. We write $H(q)=q^{k-1}J(q)$, where
$J(q)=-(j-k)q^j+jq^{j-k}-k$.
We have $J(0)<0$ and $\lim_{q\to\infty}J(q)=-\infty$.
Also, $J'(q)=-(j-k)jq^{j-1}(1-q^{-k})$. Thus, $J$ attains its maximum at $q=1$, and one has  $J(1)=0$; thus  $J(q)<0$, for $q\neq1$.
\hfill $\square$

\medskip

\noindent \it Proof of Theorem \ref{p1q1}.\rm\
From  \eqref{p1distfunc} with $q=1$, we obtain for $x\in(0,1)$,
$\lim_{n\to\infty}\mathcal{P}_n^{(1)}(\frac{\pi_1}n\le x)=x-x\log x$. This proves part (i).
From \eqref{p1dist=} with $q=1$, we have
$$
\begin{aligned}
&\mathcal{E}_n^{(1)}\pi_1=\frac1n\sum_{k=1}^nk\sum_{j=k}^n\frac1j=\frac1n\sum_{j=1}^n\frac1j\sum_{k=1}^jk=\frac1{2n}\sum_{j=1}^n(j+1)=\\
&\frac1{2n}\left(\frac12n(n+1)+n\right)=\frac{n+3}4,
\end{aligned}
$$
which proves part (ii).
Parts (iii)-(v) follow immediately from \eqref{p1dist=} with $q=1$.
\hfill $\square$

\medskip

\noindent \it Proof of Theorem \ref{p1q11n}.\rm\
We'll prove the case that $q_n=1+\frac cn$ with $c>0$; the case $c<0$ is proved similarly.
Let $x\in(0,1)$. From  \eqref{p1distfunc} we have
\begin{equation}\label{pi1asymp1cn}
\mathcal{P}_n^{(1+\frac cn)}(\frac{\pi_1}n\le x)=\frac{[xn]}n+\frac1n\left((1+\frac cn)^{[xn]}-1\right)\sum_{j=[xn]+1}^n\frac1{(1+\frac cn)^j-1}.
\end{equation}
We have $(1+\frac cn)^j\sim e^{\frac{cj}n}$, uniformly over $j\in\{[xn]+1,\cdots, n\}$. Using this with \eqref{pi1asymp1cn}, we have
\begin{equation}\label{pi1asymp1cnagain}
\lim_{n\to\infty}\mathcal{P}_n^{(1+\frac cn)}(\frac{\pi_1}n\le x)=x+(e^{cx}-1)\int_x^1\frac1{e^{cy}-1}dy.
\end{equation}
We have
\begin{equation}\label{integral}
\int_x^1\frac1{e^{cy}-1}dy=\int_x^1\frac{e^{-cy}}{1-e^{-cy}}dy=\frac1c\log(1-e^{-cy})|_x^1=\frac1c\log\frac{1-e^{-c}}{1-e^{-cx}}.
\end{equation}
Substituting \eqref{integral} in \eqref{pi1asymp1cnagain} gives
\eqref{distfunc1n}, and differentiating \eqref{distfunc1n} gives \eqref{density1n}.
We leave it to the reader to check that the density in \eqref{density1n} is monotone decreasing.
\hfill $\square$
\medskip

\noindent \it Proof of Theorem \ref{p1q1nalpha}.\rm\
Let $x\in(0,1)$.
 From  \eqref{p1distfunc} we have
\begin{equation}\label{pi1asympnalpha}
\begin{aligned}
&\mathcal{P}_n^{(1-\frac c{n^\alpha})}(\frac{\pi_1}{n^\alpha}\le x)=
\frac{[xn^\alpha]}n+\frac1n\left((1-(1-\frac c{n^\alpha})^{[xn^\alpha]}\right)
\sum_{j=[xn^\alpha]+1}^n\frac1{1-(1-\frac c{n^\alpha})^j}.
\end{aligned}
\end{equation}
Since for any $\beta\in(\alpha,1)$,  $\lim_{n\to\infty}(1-\frac c{n^\alpha})^j=0$, uniformly over $j\in[n^\beta,n]$,
we have $\lim_{n\to\infty}\frac1n\sum_{j=[xn^\alpha]+1}^n\frac1{1-(1-\frac c{n^\alpha})^j}=1$.
Using this in \eqref{pi1asympnalpha}, we obtain
$$
\lim_{n\to\infty}\mathcal{P}_n^{(1-\frac c{n^\alpha})}(\frac{\pi_1}{n^\alpha}\le x)=1-e^{-cx}.
$$
\hfill $\square$
\medskip

\noindent \it Proof of Theorem \ref{p1q}.\rm\
 From  \eqref{p1dist=} we have for  $k\in\mathbb{N}$,
\begin{equation*}\label{fixedq<1}
\lim_{n\to\infty}\mathcal{P}_n^{(q)}(\pi_1= k)=(1-q)q^{k-1}.
\end{equation*}
\hfill$\square$
\medskip

\noindent \it Proof of Theorem \ref{p1qlarge}.\rm\
We first consider the case of fixed $q>1$.
Let $x\in(0,1)$. From  \eqref{p1distfunc} we have
\begin{equation}\label{pi1asymplgq}
\mathcal{P}_n^{(q)}(\frac{\pi_1}n\le x)=\frac{[xn]}n+\frac1n(q^{[xn]}-1)\sum_{j=[xn]+1}^n\frac1{q^j-1}.
\end{equation}
For sufficiently large $n$, depending on $q$, we have
$$
\begin{aligned}
&\frac1n(q^{[xn]}-1)\sum_{j=[xn]+1}^n\frac1{q^j-1}\le \frac{2q^{[xn]}}{n}\sum_{j=[xn]+1}^\infty\frac1{q^j}=\frac2{(q-1)n}.
\end{aligned}
$$
Using this with \eqref{pi1asymplgq} gives
$$
\lim_{n\to\infty}\mathcal{P}_n^{(q)}(\frac{\pi_1}n\le x)=x,
$$
proving that $\frac{\pi_1} n$ converges in distribution to the uniform distribution on $[0,1]$.

The first equation in \eqref{cornersq} is \eqref{p1dist=}. For the second equation in \eqref{cornersq}, we have from \eqref{p1dist=},
$$
\begin{aligned}
&\mathcal{P}_n^{(q)}(\pi_1=n-k)=\frac1n(q-1)\sum_{i=0}^k\frac{q^{n-k-1}}{q^{n-k+i}-1}\sim\frac1n(q-1)\sum_{i=0}^k(\frac1q)^{i+1}=\\
&\frac1n(q-1)\frac{\frac1q-(\frac1q)^{k+2}}{1-\frac1q}=\frac1n\left(1-\frac1{q^{k+1}}\right).
\end{aligned}
$$

Now consider
the general case that $q_n>1$ is bounded and $q_n-1=\omega(\frac1n)$. Then it follows from Proposition \ref{stochdom} that for any fixed $c>1$, for sufficiently large $n$
the distribution of $\pi_1$ under $\mathcal{P}_n^{(1+\frac cn)}$ is stochastically dominated by the distribution of $\pi_1$ under $\mathcal{P}_n^{(q_n)}$. Since $\{\frac{\pi_1}n\}$ is bounded, it is tight under any sequence
of distributions $\mathcal{P}_n^{(q_n)}$.
Letting $n\to\infty$ and using Theorem \ref{p1q11n},   it follows  that any  limit point in  distribution of $\{\frac{\pi_1}n\}$ under
$\mathcal{P}_n^{(q_n)}$ stochastically dominates the distribution $F_c$ in \eqref{distfunc1n}, for any $c>0$.
The term involving $c$ on the right hand side of \eqref{distfunc1n} satisfies
$$
\frac1c(e^{cx}-1)\log\frac{1-e^{-c}}{1-e^{-cx}}\sim\frac1c(e^{cx}-1)\left(-e^{-c}+e^{-cx}\right)\stackrel{c\to\infty}{\to0}.
$$
Thus, $\lim_{c\to\infty}F_c(x)=x.$
We conclude that any limit point in distribution of
$\{\frac{\pi_1}n\}$ under
$\mathcal{P}_n^{(q_n)}$  stochastically dominates the uniform distribution on $[0,1]$.
On the other hand, choose a fixed $q>\sup_{n\in\mathbb{N}} q_n$.
By the proof  above for the case of fixed $q$, along with Proposition \ref{stochdom}, it follows that
any limit point in distribution of
$\{\frac{\pi_1}n\}$ under
$\mathcal{P}_n^{(q_n)}$ is stochastically dominated by the uniform distribution on $[0,1]$.
We thus conclude that the only limit point in distribution of
$\frac{\pi_1} n$ under $\mathcal{P}_n^{(q_n)}$ is the uniform distribution. Consequently,
$\frac{\pi_1} n$ under $\mathcal{P}_n^{(q_n)}$ converges to the uniform distribution.

\hfill
$\square$
\section{Proofs of Theorems \ref{Nkqnear1}-\ref{Nkq>1}}\label{Nkproofs}
In this section we prove the results concerning the parking function statistic $N^{(n)}_k$.
From \eqref{bijection} and  the definition of
$\mathcal{P}_n^{(q)}$, it follows that $N^{(n)}_k$ under
$\mathcal{P}_n^{(q)}$ has the distribution of
$\sum_{j=k}^n1_{\{\tilde I_{n,<j}=k\}}$,
where $\tilde I_{n,<j}$ is as in \eqref{MallowsI<k}
and $\{\tilde I_{n,<j}\}_{j=1}^n$ are independent:
\begin{equation}\label{Nkdist}
N^{(n)}_k\stackrel{\text{dist}}{=}\sum_{j=k}^n1_{\{\tilde I_{n,<j}=k\}}.
\end{equation}
Thus, the Laplace transform of $N^{(n)}_k$ under  $\mathcal{P}_n^{(q)}$ is given by
\begin{equation}\label{laplace}
\mathcal{E}_n^{(q)}e^{-tN^{(n)}_k}=\begin{cases}\prod_{j=k}^n\left(1-\frac{1-q}{1-q^j}q^{k-1}(1-e^{-t})\right),\ t\ge0; \ q\in(0,\infty)-\{1\};\\
\prod_{j=k}^n\left(1-\frac1j(1-e^{-t})\right),\ t\ge0;\ q=1.\end{cases}
\end{equation}
Using \eqref{laplace}, we now prove Theorems \ref{Nkqnear1}-\ref{Nkq>1}.
\medskip

\noindent \it Proof of Theorem \ref{Nkqnear1}.\rm\
The assumption on $q_n$  in parts (i), (ii) and (iv) of the theorem, is that there exists a constant $C>0$ such that $1-\frac Cn\le q_n\le 1+\frac Cn$.
The assumption in part (iii) is that there exists a $c\in\mathbb{R}$ such that $q_n=1+\frac cn$.
To avoid repeating similar types of analysis,
we will prove all of the parts under the assumption that $q_n=1+\frac cn$, with $c>0$.
The extension to the more general form for $q_n$ is easy to obtain from the analysis below for $q_n=1+\frac cn,\ c>0$.
From \eqref{laplace}, we have
\begin{equation}\label{laplaceNkqnear1}
\mathcal{E}_n^{(1+\frac cn)}e^{-t\frac{N^{(n)}_k}{\log n}}=\prod_{j=k}^n\left(1-\frac cn\frac1{(1+\frac cn)^j-1}(1+\frac cn)^{k-1}\left(1-e^{-\frac t{\log n}}\right)\right).
\end{equation}

We first consider part (i), where $k$ is fixed.
From the Taylor expansion, we have
$$
(1+\frac cn)^j-1=\frac{jc}n+\frac{j(j-1)}2\frac{a_{j,n}^2}{n^2},\ 0<a_{j,n}<c.
$$
Thus,
$$
\frac cn\frac1{(1+\frac cn)^j-1}=\frac c{j\left(c+\frac{j-1}{2n}a_{j,n}^2\right)}.
$$
The expression above  is uniformly bounded in $j$ and $n$. Thus,    we have
\begin{equation}\label{laplaceasymp}
\begin{aligned}
&\log\left(1-\frac cn\frac1{(1+\frac cn)^j-1}(1+\frac cn)^{k-1}\left(1-e^{-\frac t{\log n}}\right)\right)\sim\frac c{j\left(c+\frac{j-1}{2n}a_{j,n}^2\right)}\frac t{\log n},\\
& \text{uniformly over}\ j\in\{k,\cdots, n\}\ \text{as}\ n\to\infty.
\end{aligned}
\end{equation}
From \eqref{laplaceNkqnear1} and \eqref{laplaceasymp}, we have
\begin{equation}\label{loglaplaceasymp}
\log\mathcal{E}_n^{(1+\frac cn)}e^{-t\frac{N^{(n)}_k}{\log n}}\sim-\frac t{\log n}\sum_{j=k}^n\frac c{j\left(c+\frac{j-1}{2n}a_{j,n}^2\right)}.
\end{equation}
For any $\epsilon>0$, there exists a $\delta_\epsilon>0$ such that
\begin{equation}\label{epsdel}
(1-\epsilon)\frac1j\le \frac c{j\left(c+\frac{j-1}{2n}a_{j,n}^2\right)}\le \frac 1j,\ j=k,\cdots, [\delta_\epsilon n].
\end{equation}
From \eqref{epsdel},
\begin{equation}\label{liminflimsup}
\begin{aligned}
&\liminf_{n\to\infty}\frac t{\log n}\sum_{j=k}^{[\delta_\epsilon n]}\frac c{j\left(c+\frac{j-1}{2n}a_{j,n}^2\right)}\ge(1-\epsilon)t;\\
&\limsup_{n\to\infty}\frac t{\log n}\sum_{j=k}^{[\delta_\epsilon n]}\frac c{j\left(c+\frac{j-1}{2n}a_{j,n}^2\right)}\le t.
\end{aligned}
\end{equation}
We also have
\begin{equation}\label{secondpiece}
\begin{aligned}
\limsup_{n\to\infty}\frac t{\log n}\sum_{j=[\delta_\epsilon n]+1}^n\frac c{j\left(c+\frac{j-1}{2n}a_{j,n}^2\right)}=0,\ \text{for all}\ \delta_\epsilon>0.
\end{aligned}
\end{equation}
From \eqref{loglaplaceasymp}, \eqref{liminflimsup} and \eqref{secondpiece}, we conclude that $\mathcal{E}_n^{(1+\frac cn)}e^{-t\frac{N^{(n)}_k}{\log n}}=e^{-t}$,
and consequently, by the uniqueness of the Laplace transform for nonnegative distributions, part (i) follows.

We now consider  part (ii), where $k_n=\theta(n^\alpha)$. The analysis above up to and including  \eqref{epsdel} holds with $k$ replaced by $k_n$.
When we replace $k$ by $k_n$ in the  lower limit in the two sums on the left hand side of \eqref{liminflimsup}, the right hand sides get multiplied by $1-\alpha$.
This proves part (ii).

We now consider part (iii) where $k_n\sim dn$ with $d\in(0,1)$.
Similar to  \eqref{laplaceNkqnear1}, we have
$$
\mathcal{E}_n^{(1+\frac cn)}e^{-tN^{(n)}_{[dn]}}=\prod_{j=[dn]}^n\left(1-\frac cn\frac1{(1+\frac cn)^j-1}(1+\frac cn)^{[dn]-1}\left(1-e^{-t}\right)\right).
$$
Thus,
\begin{equation}\label{dncaseasymp}
\begin{aligned}
&\log \mathcal{E}_n^{(1+\frac cn)}e^{-tN^{(n)}_{[dn]}}=\sum_{j=[dn]}^n\log\left(1-\frac cn\frac1{(1+\frac cn)^j-1}(1+\frac cn)^{[dn]-1}\left(1-e^{-t}\right)\right)\sim\\
&-\frac cne^{cd}(1-e^{-t})\sum_{j=[dn]}^n\frac1{e^{\frac{cj}n}-1}\sim -ce^{cd}(1-e^{-t})\int_d^1\frac1{e^{cx}-1}dx.
\end{aligned}
\end{equation}
From \eqref{dncaseasymp} and \eqref{integral}, we have
$$
\lim_{n\to\infty}\log \mathcal{E}_n^{(1+\frac cn)}e^{-tN^{(n)}_{[dn]}}=-\left(e^{cd}\log\frac{1-e^{-c}}{1-e^{-cd}}\right)(1-e^{-t}).
$$
Thus,
$$
\lim_{n\to\infty}\mathcal{E}_n^{(1+\frac cn)}e^{-tN^{(n)}_{[dn]}}=e^{-\left(e^{cd}\log\frac{1-e^{-c}}{1-e^{-cd}}\right)(1-e^{-t})}.
$$
The right hand side above is the Laplace transform of the Poisson distribution with parameter $e^{cd}\log\frac{1-e^{-c}}{1-e^{-cd}}$; so $N^{(n)}_{[dn]}$ converges to this Poisson distribution.

We now consider part (iv) where $k_n$ satisfies $n-k_n=o(n)$. Write $k_n=n-a(n)$ with $a(n)=o(n)$. Then similar to \eqref{dncaseasymp}, we have
\begin{equation}\label{caseivasymp}
\begin{aligned}
&\log \mathcal{E}_n^{(1+\frac cn)}e^{-tN^{(n)}_{k_n}}=\sum_{j=n-a_n}^n\log\left(1-\frac cn\frac1{(1+\frac cn)^j-1}(1+\frac cn)^{n-a_n-1}\left(1-e^{-t}\right)\right)\sim\\
&-\frac cne^c(1-e^{-t})\sum_{j=n-a_n}^n\frac1{e^{\frac{cj}n}-1}.
%-ce^{cd}(1-e^{-t})\int_d^1\frac1{e^{cx}-1}dx.
\end{aligned}
\end{equation}
From \eqref{caseivasymp}, we conclude that
$\lim_{n\to\infty}\log \mathcal{E}_n^{(1+\frac cn)}e^{-tN^{(n)}_{k_n}}=0$, and thus,\newline
 $\lim_{n\to\infty}\mathcal{E}_n^{(1+\frac cn)}e^{-tN^{(n)}_{k_n}}=1$. Therefore, $N^{(n)}_{k_n}$ converges to zero in distribution.
 \hfill $\square$

\medskip

\noindent \it Proof of Theorem \ref{Nkq<1}.\rm\ We first consider part (i) where $k$ is fixed. From \eqref{laplace}, we have
\begin{equation}\label{q<1Nk}
\begin{aligned}
&\log\mathcal{E}_n^{(q)}e^{-t\frac{N^{(n)}_k}n}=\sum_{j=k}^n\log\left(1-\frac{1-q}{1-q^j}q^{k-1}(1-e^{-\frac tn})\right)\sim\\
&-(1-q)q^{k-1}\frac tn\sum_{j=k}^n\frac1{1-q^j}.
\end{aligned}
\end{equation}
Thus,
$\lim_{n\to\infty}\log\mathcal{E}_n^{(q)}e^{-t\frac{N^{(n)}_k}n}=-(1-q)q^{k-1}t$, and
$$
\lim_{n\to\infty}\mathcal{E}_n^{(q)}e^{-t\frac{N^{(n)}_k}n}=e^{-(1-q)q^{k-1}t},
$$
which proves that $\frac{N^{(n)}_k}n$ converges in distribution to $(1-q)q^{k-1}$.

We now consider part (ii) where $\lim_{n\to\infty}nq^{k_n}=\infty$. Similar to \eqref{q<1Nk},
we have
\begin{equation}\label{q<1Nkagain}
\begin{aligned}
&\log\mathcal{E}_n^{(q)}e^{-t\frac{N^{(n)}_{k_n}}{nq^{k_n}}}=\sum_{j=k_n}^n\log\left(1-\frac{1-q}{1-q^j}q^{k_n-1}(1-e^{-\frac t{nq^{k_n}}})\right)\sim\\
&-(1-q)q^{k_n-1}\frac t{nq^{k_n}}\sum_{k=k_n}^n\frac1{1-q^j}.
\end{aligned}
\end{equation}
Since $k_n=o(n)$, it  follows from \eqref{q<1Nkagain}  that $\lim_{n\to\infty}\log\mathcal{E}_n^{(q)}e^{-t\frac{N^{(n)}_{k_n}}{nq^{k_n}}}=\frac{1-q}qt$, and
$$
\lim_{n\to\infty}\mathcal{E}_n^{(q)}e^{-t\frac{N^{(n)}_{k_n}}{nq^{k_n}}}=e^{-\frac{1-q}qt},
$$
which proves that
$\frac{N^{(n)}_{k_n}}{nq^{k_n}}$ converges in distribution to $\frac{1-q}q$.

We now consider parts (iii) and (iv) together. Thus, we assume that $\lim_{n\to\infty}nq^{k_n}=L\in[0,\infty)$.
Then similar to \eqref{q<1Nk} and \eqref{q<1Nkagain}, we have
\begin{equation}\label{q<1Nkn}
\begin{aligned}
&\log\mathcal{E}_n^{(q)}e^{-tN^{(n)}_{k_n}}=\sum_{j=k_n}^n\log\left(1-\frac{1-q}{1-q^j}q^{k_n-1}(1-e^{-t})\right)\sim\\
&-(1-q)q^{k_n-1}(1-e^{-t})\sum_{j=k_n}^n\frac1{1-q^j}=-(1-q)nq^{k_n-1}(1-e^{-t})\frac1n\sum_{j=k_n}^n\frac1{1-q^j}.
\end{aligned}
\end{equation}
Since $k_n=o(n)$, it  follows  from \eqref{q<1Nkn} that
$\lim_{n\to\infty}\log\mathcal{E}_n^{(q)}e^{-tN^{(n)}_{k_n}}=\frac{1-q}qL(1-e^{-t})$,
and
\begin{equation}\label{q<1Nknagain}
\lim_{n\to\infty}\mathcal{E}_n^{(q)}e^{-tN^{(n)}_{k_n}}=e^{-\frac{1-q}qL(1-e^{-t})}.
\end{equation}
If $L>0$, the right hand of \eqref{q<1Nknagain} is the Laplace transform of the Poisson distribution with parameter $\frac{1-q}qL$;
so $N^{(n)}_{k_n}$ converges to this Poisson distribution.
If $L=0$, then it follows from \eqref{q<1Nknagain} that $N^{(n)}_{k_n}$ converges to zero in distribution.

\medskip

\noindent \it Proof of Theorem \ref{Nkq>1}.\rm\
We begin with part (i).
By \eqref{Nkdist} and \eqref{MallowsI<k}, it follows that $N^{(n)}_k$ converges in distribution to the distribution of
$\sum_{j=k}^\infty Z_{j;k}$, as in the statement of the theorem.
Since $E\sum_{j=k}^\infty Z_{j:k}=\sum_{j=k}^\infty \frac{(q-1)q^{k-1}}{q^j-1}<\infty$, it follows that
the support of $\sum_{j=k}^\infty Z_{j;k}$ is $\mathbb{Z}^+$.

We now turn to part (ii). By \eqref{Nkdist} and \eqref{MallowsI<k}, we can write  for any $M\in\mathbb{N}$ and sufficiently large $n$,
\begin{equation*}\label{breakintwo}
N^{(n)}_k\stackrel{\text{dist}}{=}\sum_{i=0}^M Z_{k_n+i;k_n}+\sum_{i=M+1}^{n-k_n}Z_{k_n+i;k_n}:=I_M+II_{M,n}.
\end{equation*}
From the definition of $Z_{k_n+i;k_n}$ and $Y_i$, clearly, $I_M$ converges in distribution as $n\to\infty$ to $\sum_{i=0}^M Y_i$,
where the random variables  $\{Y_i\}_{i=0}^\infty$ are as in the statement of the theorem.
Also,
$$
EII_{M,n}=\sum_{i=M+1}^{n-k_n}\frac{(q-1)q^{k_n-1}}{q^{k_n+i}-1}.
$$
Thus, for sufficiently large $n$,
$$
EII_{M,n}\le 2(q-1)\sum_{i=M+1}^\infty q^{-1-i}=2q^{-M}.
$$
Thus, $\lim_{M\to\infty}\limsup_{n\to\infty}EII_{M,n}=0$.
This completes the proof of part (ii).

We now turn to part (iii). By \eqref{Nkdist} and \eqref{MallowsI<k}, we can write
$$
N^{(n)}_{n-k}\stackrel{\text{dist}}{=}\sum_{i=0}^k Z_{n-i;n-k}.
$$
Thus, from the definition of $Z_{n-i;n-k}$ and $Y_{k-i}$, clearly $N^{(n)}_{n-k}$ converges in distribution as $n\to\infty$ to $\sum_{i=0}^k Y_i$.

For part (iv), we write
$\sum_{j=k}^\infty Z_{j;k}=\sum_{i=0}^\infty Z_{k+i;k}$ and note that \newline $\lim_{k\to\infty}Z_{k+i;k}\stackrel{\text{dist}}{=}Y_i$.
\hfill $\square$

\end{document}